\documentclass[review]{elsarticle}
\usepackage{lineno}
\usepackage{bbold}
\usepackage{verbatim}
\usepackage{url}
\usepackage{graphicx}
\usepackage{float}
\usepackage{mathrsfs}
\usepackage{epstopdf}
\usepackage[all,color]{xy}
\usepackage{cancel}
\usepackage{MnSymbol}
\usepackage{wasysym, stmaryrd}
\usepackage{tikz}
\usepackage[margin=1in]{geometry}

\usetikzlibrary{arrows}
\modulolinenumbers[5]
\allowdisplaybreaks
\newtheorem{theorem}{Theorem}

\newtheorem{corollary}[theorem]{Corollary}
\newtheorem{example}[theorem]{Example}
\newtheorem{lemma}[theorem]{Lemma}

\newproof{proof}{Proof}
\numberwithin{equation}{subsection}
\numberwithin{theorem}{subsection}

\newcommand{\ignore}[1]{}
\newcommand*{\vcenteredhbox}[1]{\begin{tabular}{@{}c@{}}#1\end{tabular}}

\begin{document}

\begin{frontmatter}

\title{The Determinant of $\{\pm 1\}$-Matrices and Oriented Hypergraphs}

\author[add2]{Lucas J. Rusnak\corref{mycorrespondingauthor}}\ead{Lucas.Rusnak@txstate.edu}
\author[add2,add3]{Josephine Reynes}
\author[add1]{Russell Li}
\author[add1,add4]{Eric Yan}
\author[add1]{Justin Yu}

\address[add2]{Department of Mathematics, Texas State University, San Marcos, TX 78666, USA}

\address[add3]{Combinatorics and Optimization, University of Waterloo, Waterloo ON, N2L 3G1, Canada}

\address[add4]{Harvard College, Cambridge, MA 02138, USA}

\address[add1]{Mathworks, Texas State University, San Marcos, TX 78666, USA}

\cortext[mycorrespondingauthor]{Corresponding author}

\begin{abstract}
The determinants of $\{\pm 1\}$-matrices are calculated by via the oriented hypergraphic Laplacian and summing over an incidence generalization of vertex cycle-covers. These cycle-covers are signed and partitioned into families based on their hyperedge containment. Every non-edge-monic family is shown to contribute a net value of $0$ to the Laplacian, while each edge-monic family is shown to sum to the absolute value of the determinant of the original incidence matrix. Simple symmetries are identified as well as their relationship to Hadamard's maximum determinant problem. Finally, the entries of the incidence matrix are reclaimed using only the signs of an adjacency-minimal set of cycle-covers from an edge-monic family.

\end{abstract}


\end{frontmatter}

%%%%%%%%

\section{Introduction and Background}

An oriented hypergraph \cite{AH1} is a generalization of bidirected graphs and signed graphs \cite{MR0267898, OSG} that allows for integer incidence matrices to be studied locally using signed graphs. Oriented hypergraphs allow for the extension of graphic concepts to integer matrices this includes spectral properties \cite{Reff6,Reff7,Mulas3,Mulas1,Reff2}, applications to chemical hypergraphs \cite{Mulas4,Mulas2}, and Hadamard matrices \cite{Reff3,Reff5}. We examine the connection between the determinant of $\{\pm 1\}$-matrices as uniform oriented hypergraphs and the generalized cycle-cover formula for integer matrix characteristic polynomials from \cite{OHSachs,OHMTT} to determine a new, highly symmetric, combinatorial interpretation of matrix determinants and discuss its relationship to Hadamard's maximum determinant problem \cite{Had93}.

\begin{figure}[H]
    \centering
    \includegraphics{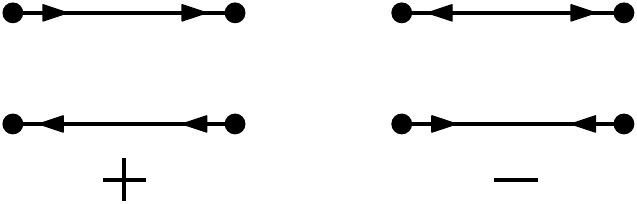}
    \caption{Bidirected orientations of positive and negative edges.}
    \label{fig:BiDirEdge}
\end{figure}

A bidirected edge \cite{MR0267898} has an orientation at each incidence and represents an orientation of a signed graphic edge \cite{OSG}. Positive edges are oriented to ``travel through'' the edge, while negative edges are oriented so the incidences oppose each other, as depicted in Figure \ref{fig:BiDirEdge}.

\begin{figure}[H]
    \begin{center}
\vcenteredhbox{\includegraphics{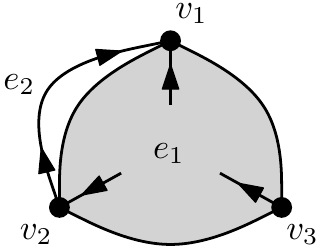}} \qquad  \vcenteredhbox{$\mathbf{H}_G=
\left[ 
\begin{array}[c]{cc}
1 & 1 \\
1 & -1 \\
-1 & 0 
\end{array}
\right]$}
\end{center}
    \caption{An oriented hypergraph and its incidence matrix.}
    \label{fig:OH}
\end{figure}

An oriented hypergraph extends the notion of bidirectedness to multidirectedness for hyperedges, where the incidences are oriented, and adjacency signs are determined locally as though they were a signed graph. For the oriented hypergraph in Figure \ref{fig:OH} (left), observe that the adjacency from $v_1$ to $v_2$ along edge $e_1$ is negative, while the adjacency from $v_1$ to $v_2$ along edge $e_2$ is positive.

The \emph{incidence matrix} of an oriented hypergraph $G$ is the $V \times E$ matrix $\mathbf{H}_{G}$ where the $(v,e)$-entry is the sum of the orientation values of the incidences between $v$ and $e$ ($+1$ for entrant and $-1$ for salient). Figure \ref{fig:OH} (right) depicts the incidence matrix of the given oriented hypergraph. 

Subsection \ref{ssec:OHIntro} contains an abbreviated background on oriented hypergraphs and signed adjacency, we refer the reader to \cite{AH1,OH1} for additional background on the algebraic and structural basics of oriented hypergraphs. Section \ref{sec:OHLaplacian} reviews how the determinant of any integer matrix Laplacian is calculated based on the characteristic polynomial coefficient theorem and the notion of generalized cycle-covers called \emph{contributors} from \cite{OHSachs,OHMTT}. In Subsection \ref{ssec:ReduceToEM} we prove that the determinant of the Laplacian of any $\{\pm 1\}$-matrix is determined only by the edge-monic contributors which contain a representative from each hyperedge. Subsection \ref{ssec:EMIdentAndHeadTail} addresses the entangled nature of contributor calculations used for the determinant of a Laplacian.

Section \ref{sec:DetOfH} begins with an example summarizing the results from Section \ref{sec:OHLaplacian} by doing a partial calculation for the Laplacian determinant by examining a single edge-monic class of contributors. This class is observed to sum (up to sign) to the determinant of the original incidence matrix, and is proved in general in Subsection \ref{ssec:OneEM}. Simple formulas for determining the determinant of any $\{\pm 1\}$-matrix based on the contributor signs are produced for a single edge-monic contributor class as well as the set of edge-monic contributor classes.

Section \ref{sec:HadAndMax} acknowledges the connections to Hadamard's maximum determinant problem \cite{Had93} by varying the oriented hypergraphic orientations to minimizing contributors of a given sign. Subsection \ref{ssec:01Equiv} demonstrates that the equivalence of determinants between $\{\pm 1\}$-matrices and $\{0,+1\}$-matrices is related to negative digon placement in a specific class of fundamental circles in the associated oriented hypergraph. In Subsection \ref{ssec:ReclaimH} we conclude by showing that the contributor signs corresponding to an adjacency-minimal collection of fundamental circles of a single edge-monic class allows for the the unique reconstruction of the original $\{\pm 1\}$-matrix, assuming the rows and columns have been negated appropriately so that the first row and column consist of only $+1$'s.

\subsection{Oriented Hypergraphs}
\label{ssec:OHIntro}

Formally, an \emph{oriented hypergraph} is a sextuple $G=(V, E, I, \varsigma, \omega, \sigma)$ consisting of a set of vertices $V$, a set of edges $E$, a set of incidences $I$, two incidence end maps $\varsigma:I\to V$, and $\omega:I\to E$, and an incidence orientation function $\sigma:I\rightarrow\{+1,-1\}$. Incidence orientations of $+1$ are indicated by an arrow entering a vertex, and $-1$ exiting a vertex. The quintuple $(V, E, I, \varsigma, \omega)$ is the underlying \emph{incidence hypergraph} of $G$.

\begin{figure}[H]
    \centering
    \includegraphics{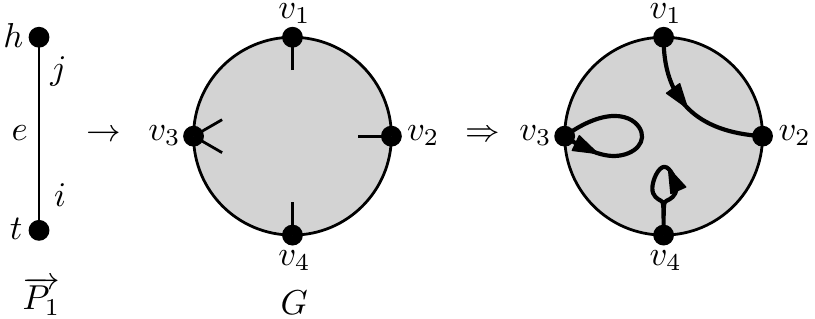}
    \caption{A directed adjacency between $v_1$ and $v_2$, a loop at $v_3$, and a backstep at $v_4$.}
    \label{fig:ABL}
\end{figure}

Let the sequence $\overrightarrow{P}_{1}= (t,i,e,j,h)$ denote a path of length $1$ from tail-vertex $t$ to head-vertex $h$ with edge $e$ and incidences $i$ and $j$ such that $i$ is between $t$ and $e$, and $j$ is between $h$ and $e$. The incidences $i$ and $j$ are the \emph{tail-incidence} and \emph{head-incidence}, respectively. It was shown in \cite{AH1} that the common graphic matrices of the incidence, adjacency, degree, and Laplacian extend to all integer matrices through its oriented hypergraphic representation via incidence preserving homomorphisms $p:\overrightarrow{P}_{1} \to G$. A \emph{directed adjacency} is an incidence-monic image of $p$, a \emph{loop} is directed adjacency that is not vertex-monic, and a \emph{backstep} is a non-incidence-monic image of $p$. Here, \emph{monic} means the object appears at most one time, see Figure \ref{fig:ABL} for an example of a directed adjacency, loop, and backstep.

Observe that the incidence interpretation is able to discern the orientations of loops, which canonically produce a value of $2$ on the diagonal of the adjacency matrix. However, a backstep travels from a vertex into an edge and immediately returns to the original vertex along the same incidence, the repeated incidence prevents determining an orientation.

To avoid confusion with cycle-spaces and matroidal circuits, the image of a simple graphic cycle that does not repeat a vertex or an edge is called a \emph{circle}. A \emph{digon} is a circle containing exactly two adjacencies. We reserve the term ``cycle'' for a permutation cycle, which manifests as a closed sequence of $p$ maps into $G$, and need not necessarily be a circle in the oriented hypergraph. The \emph{sign of an adjacency (backstep) in $G$} is $-\sigma(p(i))\sigma(p(j))$, and the sign of any sequence of adjacencies is the product of the individual adjacency signs --- as these locally resemble bidirected edges we refer to the ``locally signed graphic structure.''

\section{Oriented Hypergraphs, Laplacians and Determinants}
\label{sec:OHLaplacian}

The \emph{Laplacian matrix} of an oriented hypergraph $G$ is $\mathbf{L}_{G}:=\mathbf{H}_{G} \mathbf{H}_{G}^{T}=\mathbf{D}_{G}-\mathbf{A}_{G}$, where the latter equality uses the incidence-degree and signed adjacency matrix, see \cite{AH1}. The determinant of the oriented hypergraphic Laplacian can be calculated using an incidence-based generalization of cycle-covers introduced in \cite{OHSachs}, called contributors. Formally, a \emph{contributor of $G$} is a $\left| V \right|$-tuple of incidence preserving path maps from a disjoint union of $\overrightarrow{P}_{1}$'s defined by $c:\dcoprod \limits_{v\in V}\overrightarrow{P}_{1}\rightarrow G$ such that $c(t_{v})=v$ and $\{c(h_{v})\mid v\in V\}=V$. The $\overrightarrow{P}_{1}$ images simply cover the vertex set twice, once with a tail-vertex and once with a head-vertex. A \emph{component} of a contributor is a restriction of the contributor map that corresponds to a single permutation cycle. Backsteps and loops are representations of $1$-cycles.

Contributors are built to mirror permutations from $V \to V$, however, far more are produced due to the incidence structure. They represent the finest (non-redundant) set of integer sums needed to determine integer matrix immanants, hence they are the fundamental calculations that ``contribute'' to the determinant and permanent calculations. The determinant of the oriented hypergraphic Laplacian is calculated according the to the following theorem:

\begin{theorem}[\cite{OHSachs}]
\label{OHSachs1}
Let $G$ be an oriented hypergraph with Laplacian matrix $\mathbf{L}_{G}$, then
\begin{align*}
    \det(\mathbf{L}_{G}) = \dsum\limits_{c\in {\mathfrak{C}}(G)}(-1)^{pc(c)},
\end{align*}
where $\mathfrak{C}(G)$ is the set of contributors of $G$, and $pc(c)$ is the number of positive components in contributor $c$.
\end{theorem}

\begin{figure}[H]
    \centering
    \includegraphics{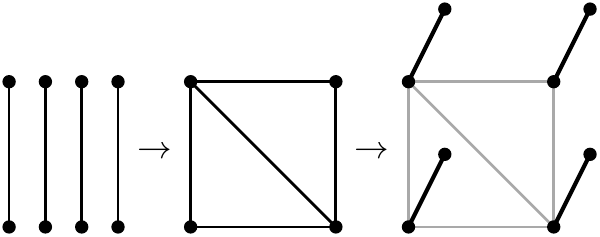}
    \caption{Path images of tail-vertices covering the vertex set of $G$.}
    \label{fig:Cont1}
\end{figure}

\begin{figure}[H]
    \centering
    \includegraphics{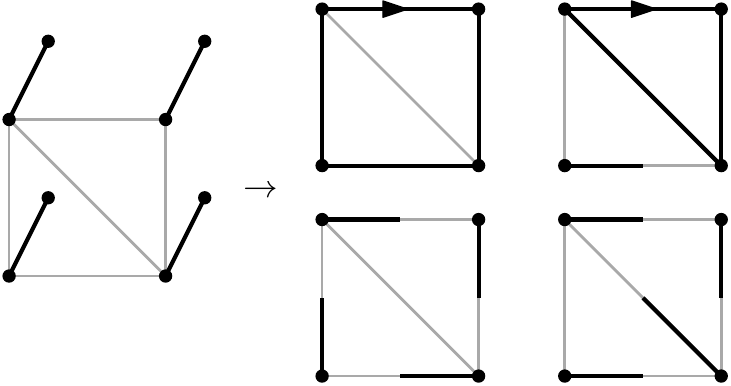}
    \caption{Four contributors formed by forcing the head-vertices to cover the vertex set. Backsteps appear as half-edges.}
    \label{fig:Cont2}
\end{figure}

An example of building contributors is illustrated on a graph in Figures \ref{fig:Cont1} and \ref{fig:Cont2}. Take a copy of $\overrightarrow{P}_{1}$ for each vertex and map each tail to cover the vertices of $G$, as shown in Figure \ref{fig:Cont1}.

To complete the contributor, map the remaining parts of each $\overrightarrow{P}_{1}$ into an incident edge then to a vertex to again cover the vertices of $G$. Four examples with direction of maps appear in Figure \ref{fig:Cont2}. Both bottom contributors contain only backsteps and are different manifestations of the identity permutation. Also, observe that each of the identity-contributors in Figure \ref{fig:Cont2} shares the image of the tail-incidence with the contributor directly above it. 

It was shown in \cite{OHMTT} that if $G$ consists of all $2$-edges, the tail-incidence images induce a partition on all contributors into Boolean lattices ordered by cycle containment, and this produces a re-interpretation of the All-Minors Matrix Tree Theorem from \cite{Seth1} using the determinant formula from Theorem \ref{OHSachs1}. We introduce an extension of this partition that, when applied to oriented hypergraphs, provides a combinatorial characterization of both the Laplacian determinant as well as the original incidence matrix. 

Two contributors are said to be \emph{tail-equivalent} if the image of their tail-incidences agree. The elements of a tail-equivalence class are partially ordered by $c \leq c'$ if (1) the set of circles of $c$ are contained in the set of circles of $c'$, or (2) the set of incidences are equal and $c$ has more components than $c'$. Clearly, there is exactly one minimal contributor in each tail-equivalency class, it consists of only backstep maps, and is associated with the identity permutation. This ordering provides combinatorial context for contributors in terms of the signless Stirling numbers of the first kind, $s(n,k)$. For example, the contributors of a single incidence-simple $n$-edge have $s(n,k)$ elements of rank $n-k$. \emph{Head-equivalence} is defined analogously based on head-incidence agreement.

\subsection{Laplacians of $\{\pm{1}\}$-matrices}
\label{ssec:ReduceToEM}

\begin{figure}[H]
    \centering
    \includegraphics{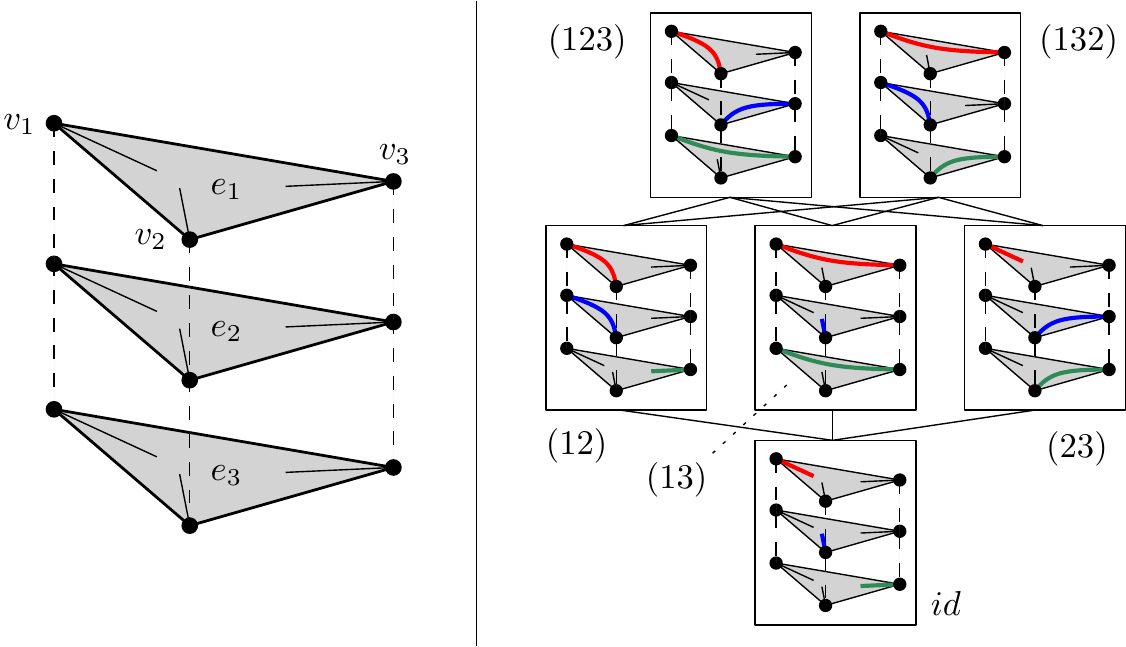}
    \caption{Left: A $3$-full incidence hypergraph; Right: A tail-equivalence class.}
    \label{fig:BaseCont}
\end{figure}

An incidence hypergraph with $\lvert V \rvert = \lvert E \rvert = n$, that is also $n$-regular and $n$-uniform is called \emph{$n$-full}, and we treat $\{\pm 1\}$-matrices as incidence matrices of $n$-full oriented hypergraphs. As drawing $n$-full hypergraphs is problematic, they will be depicted with their edges as a ``stack'' where each vertex appears in each edge, and the identification of these vertices produces the oriented hypergraph (see \cite{Reff3, Reff5}). Figure \ref{fig:BaseCont} depicts a $3$-full oriented hypergraph (left), with one tail-equivalence class (right).

For an $n$-full hypergraph there are $n!$ elements in each tail-equivalence class, one for each permutation, as each $\overrightarrow{P}_{1}$ has the option of continuing to any vertex along the edge determined by its tail-incidence. There are also $n^n$ different tail-equivalence classes, one for each contributor that consists of only backsteps. Thus, there are $n^n n!$ total contributors to determine the determinant of the Laplacian of an $n$-full oriented hypergraph using Theorem \ref{OHSachs1}.

We now partition tail-equivalence classes as follows: (1) the \emph{edge-monic} case where the image of each tail-incidence resides in a different edge, and (2) the non-edge-monic case where two tail-incidences reside in a common edge. Observe that the identity-contributor in Figure \ref{fig:BaseCont} is edge-monic as there is a backstep in each edge. Let $\mathfrak{A}$ denote the set of tail-equivalence classes, let $\mathfrak{A}_1$ be the set of edge-monic tail-equivalence classes, and let $\mathfrak{A}_2$ be the set of non-edge-monic tail-equivalence classes. We now show that the non-edge-monic tail-equivalence classes must sum to $0$ over any orientation of an $n$-full incidence hypergraph. 

\begin{theorem}
For an $n$-full oriented hypergraph $G$, the sum of the non-edge-monic contributors is zero regardless of orientation, and
\begin{align*}
    \det(\mathbf{L}_{G}) = \dsum\limits_{ \mathcal{A} \in \mathfrak{A}_1}\dsum\limits_{c \in \mathcal{A}} (-1)^{pc(c)}.
\end{align*}
\label{t:NonMonicIsZero}
\end{theorem}
\begin{proof}
From Theorem \ref{OHSachs1} we sum the contributors over their activation classes:
\begin{align*}
    \det(\mathbf{L}_{G}) &= \dsum\limits_{c\in {\mathfrak{C}}
(G)}(-1)^{pc(c)} = \dsum\limits_{ \mathcal{A} \in \mathfrak{A}}\dsum\limits_{c \in \mathcal{A}} (-1)^{pc(c)} \\ 
&= \dsum\limits_{ \mathcal{A} \in \mathfrak{A}_1}\dsum\limits_{c \in \mathcal{A}} (-1)^{pc(c)} + \dsum\limits_{ \mathcal{A} \in \mathfrak{A}_2}\dsum\limits_{c \in \mathcal{A}} (-1)^{pc(c)},
\end{align*}
and we show that the second sum is zero. Specifically, for any $\mathcal{A} \in \mathfrak{A}_2$
\begin{align*}
\dsum\limits_{c \in \mathcal{A}} (-1)^{pc(c)} = 0.
\end{align*}

Let $\mathcal{A} \in \mathfrak{A}_2$, since $\mathcal{A}$ is a non-edge-monic tail-equivalence class, there exist distinct vertices $v$ and $w$ such that $c(e_v)=c(e_w)$ for every contributor $c$. Let the permutation $\pi'=\pi\cdot(vw)$ be the permutation obtained by multiplying by the transposition $(vw)$; clearly, $(\pi')'=\pi$. Also, since each tail-equivalence class contains every permutation, let $c_{\pi}$ be the contributor of $\mathcal{A}$ corresponding to permutation $\pi$. 
    
Since $\mathcal{A}$ is a tail-equivalence class the set of tail-incidences of $c_{\pi}$ and $c_{\pi'}$ are identical. Let the head maps $c_{\pi}(h_{v})=a$ and $c_{\pi}(h_{w})=b$, then we see that $c_{\pi'}(h_{v})=b$ and $c_{\pi'}(h_{w})=a$. Moreover, their respective head-incidences of $v$ and $w$ belong to the same edge so $c_{\pi}$ and $c_{\pi'}$ also have the same head-incidences. Thus, the incidences of $c_{\pi}$ and $c_{\pi'}$ are identical.
    
\textit{Case 1:} One of $c_{\pi}$ or $c_{\pi'}$ has $v$ and $w$ in a disjoint $2$-cycle using the same adjacency twice, and the other has a backstep at $v$ and $w$ within the same edge. In this case $c_{\pi}$ and $c_{\pi'}$ differ by a single positive circle. 
    
\textit{Case 2:} Neither $c_{\pi}$ and $c_{\pi'}$ have $v$ and $w$ in a disjoint $2$-cycle using the same adjacency twice. In this case $c_{\pi}$ and $c_{\pi'}$ have the same set of incidences but differ by a single adjacency. Thus, they differ by a sign.
    
Summing over all the contributors of $\mathcal{A}$ we get
    \begin{align*}
        2\sum_{c\in\mathcal{A}}(-1)^{pc(c)} = \sum_{\pi}\left[(-1)^{pc(c_\pi)}+(-1)^{pc(c_{\pi'})}\right] = 0,
    \end{align*}
and 
\begin{align*}
\dsum\limits_{ \mathcal{A} \in \mathfrak{A}_1}\dsum\limits_{c \in \mathcal{A}} (-1)^{pc(c)} + \dsum\limits_{ \mathcal{A} \in \mathfrak{A}_2}\dsum\limits_{c \in \mathcal{A}} (-1)^{pc(c)} = \dsum\limits_{ \mathcal{A} \in \mathfrak{A}_1}\dsum\limits_{c \in \mathcal{A}} (-1)^{pc(c)}.
\end{align*} \qed
\end{proof}

\subsection{Edge-monic Identifiers and Head-tail Symmetry}
\label{ssec:EMIdentAndHeadTail}

From Theorem \ref{t:NonMonicIsZero} we know that the determinant of the Laplacian of $n$-full oriented hypergraphs is determined by the sum of only those contributors in the edge-monic tail-equivalence classes. While we rely on the symmetry of the tail-equivalence classes for $n$-full oriented hypergraphs, this raises the question if the non-edge-monic tail-equivalence classes vanish in general. We now examine the symmetry between edge-monic tail-equivalence classes and their cycle reversals, edge-monic head-equivalence classes. As previously discussed, for an $n$-full hypergraph there are $n!$ edge-monic tail-equivalence classes, each containing $n!$ contributors. However, these edge-monic tail-equivalence classes are not identical, as they each use different sets of incidences and the contributor signs of $(-1)^{pc(c)}$ in Theorem \ref{OHSachs1} need to calculated individually, except for those that use the exact same adjacencies regardless of direction used, these contributors are said to be \emph{adjacency equivalent}.

\begin{figure}[H]
    \centering
    \includegraphics{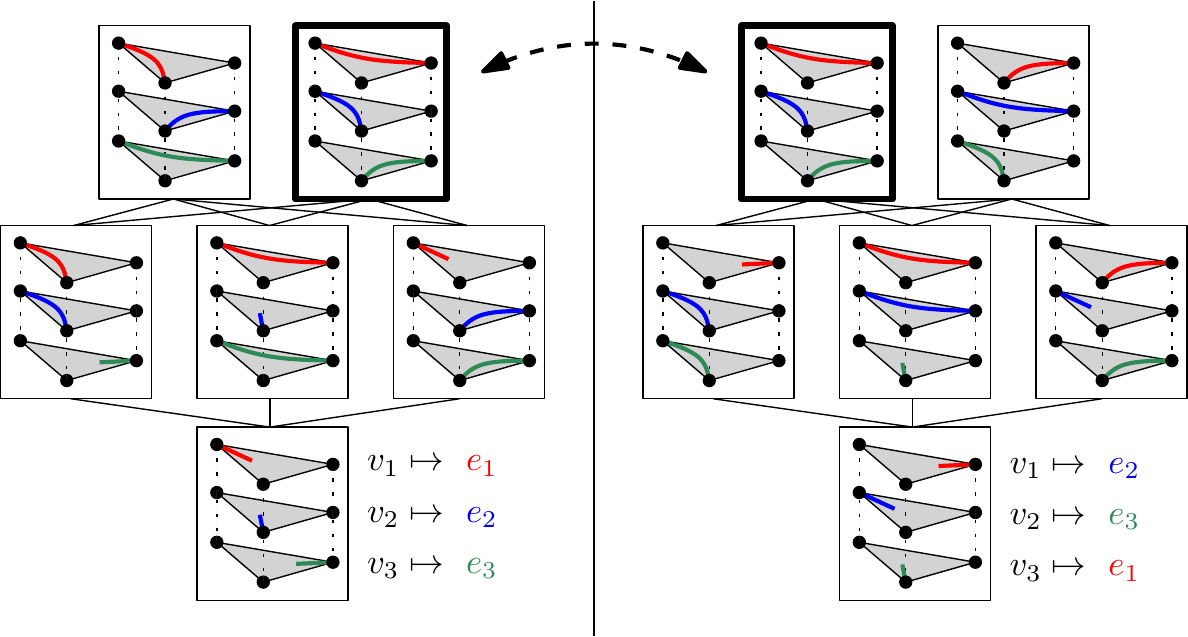}
    \caption{Two tail-equivalence classes $\mathcal{A}_{id}$ (left) and $\mathcal{A}_{(123)}$ (right) as determined by $v_k \mapsto e_{\ell}$.}
    \label{fig:EM}
\end{figure}

First, to distinguish between different edge-monic tail-equivalence classes of an $n$-full hypergraph observe that in each class the images of the tail-incidences determine ends $v_k$ and $e_{\ell}$ that induce a permutation on the indices $k \mapsto \ell$ for $k, \ell \in \{1,2,\ldots,n\}$. Call this intermediary permutation the \emph{identifier} of $\mathcal{A}$ in the set of edge-monic tail-equivalence classes $\mathfrak{A}_1$, and let $\mathcal{A}_{\alpha}$ denote the edge-monic tail-equivalence class with identifier $\alpha$. Clearly, this is a direct result of being $n$-full and being an edge-monic tail-equivalence class; $n$-fullness ensures the same number of vertices and edges, while edge-monicness ensures injectivity. Figure \ref{fig:EM} depicts two edge-monic tail-equivalence classes on the $3$-full incidence hypergraph from Figure \ref{fig:BaseCont} with identifiers $id$ and $(123)$. In Figure \ref{fig:EM} there is also a pair of adjacency-equivalent contributors, the contributors $c_{(132)} \in \mathcal{A}_{id}$ and $c_{(123)} \in \mathcal{A}_{(123)}$ (darker boxed), albeit in the reverse direction. Two adjacency-equivalent contributors that are direction reversals are said to be \emph{adjacency-inverses} as they necessarily correspond to inverse permutations.

Consider two different identifiers $\alpha$ and $\beta$ on the same $n$-full incidence hypergraph. There is exactly one such element in each of $\mathcal{A}_{\alpha}$ and $\mathcal{A}_{\beta}$ that are adjacency-inverses. These are determined by the composition 
\begin{align*}
    \xymatrix{V\ar[r]^{\alpha} & E\ar[r]^{\beta^{-1}} & V},
\end{align*}
and its inverse. These are the precisely the contributor in $\mathcal{A}_{\alpha}$ with tails determined by $\alpha$ and heads determined by $\beta$, and the contributor in $\mathcal{A}_{\beta}$ with tails determined by $\beta$ and heads determined by $\alpha$. 

\begin{lemma}
    Given any $\mathcal{A}_{\alpha}$ and $\mathcal{A}_{\beta}$ in $\mathfrak{A}_1$ in an $n$-full incidence hypergraph, there is a exactly one adjacency-equivalent pair of contributors, $c_{ \alpha^{-1}\beta} \in \mathcal{A}_{\beta}$ and $c_{\beta^{-1} \alpha} \in \mathcal{A}_{\alpha}$. Moreover, these elements are adjacency-inverses.
\end{lemma}
\begin{proof}
This is trivial, for if there were two such contributors in $\mathcal{A}_{\alpha}$ with the same head-incidences determined by $\beta$, they necessarily have to be the same contributor as their tail-incidences already coincide. A similar argument for $\mathcal{A}_{\beta}$ shows there is only one contributor with tails determined by $\beta$ and heads determined by $\alpha$. The fact they are adjacency-inverses is immediate from the fact that $\beta^{-1} \alpha$ and $\alpha^{-1}\beta$ are inverses. \qed
\end{proof}

The edge-monic tail-equivalence classes of $n$-full incidence hypergraphs are highly intertwined. Using adjacency-equivalence we are able to determine the signs of contributors in other tail-equivalence classes. However, while these contributors that belong to different tail-equivalence classes, together they form a head-equivalence class.

\begin{corollary}
    Let $\mathcal{A}_{\alpha} \in \mathfrak{A}_1$ be a tail-equivalence class in an $n$-full incidence hypergraph. The set of contributors that are adjacency-inverses to the elements of $\mathcal{A}_{\alpha}$, along with $c_{id} \in \mathcal{A}_{\alpha}$, is an edge-monic head-equivalence class. Moreover, all edge-monic head-equivalence classes are formed via adjacency-inverses of edge-monic tail-equivalence classes.
    \label{cor:Transversal}
\end{corollary}

\begin{figure}[H]
    \centering
    \includegraphics{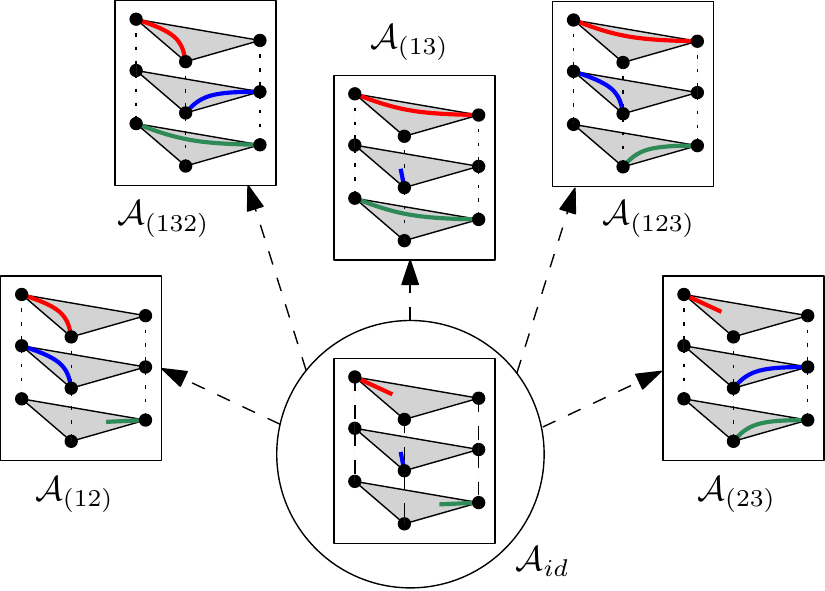}
    \caption{The elements of $\mathcal{A}_{id}$ when reversed produce an edge-monic head-equivalence class as well as a transversal of contributors from each edge-monic tail-equivalence class.}
    \label{fig:PoP2}
\end{figure}

The symmetry between edge-monic tail-equivalence classes and edge-monic head-equivalence classes raises the question is there a ``best'' class to work with when calculating the determinant of the Laplacian of an $n$-full oriented hypergraph?

\section{Determinants of $\{\pm{1}\}$-matrices}
\label{sec:DetOfH}

\begin{figure}[H]
    \begin{center}
\vcenteredhbox{$\mathbf{H}_G=
\left[ 
\begin{array}[c]{ccc}
1 & 1 & 1\\
1 & -1 & 1\\
1 & 1 & -1
\end{array}
\right]$} \qquad 
\vcenteredhbox{\includegraphics{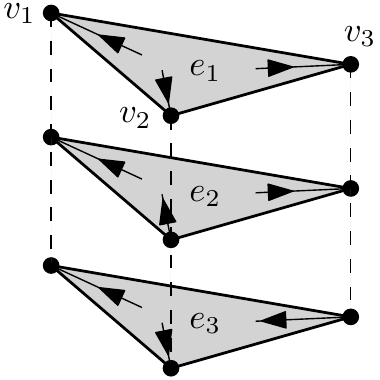}}  
\end{center}
    \caption{A $3 \times 3$ $\{\pm{1}\}$-matrix and its associated $3$-full oriented hypergraph.}
    \label{fig:OHExA}
\end{figure}

Theorem \ref{t:NonMonicIsZero} implies we only need to examine edge-monic tail-equivalence classes. We preface this section with an example of calculating the signed sum of contributors in a single  edge-monic tail-equivalence class. Consider the $3 \times 3$ matrix and its associated oriented hypergraph in Figure \ref{fig:OHExA}, and the edge-monic tail-equivalence class $\mathcal{A}_{id}$ from Figure \ref{fig:BaseCont}.

The signs of the contributors in Figure \ref{fig:BaseCont}, as determined by $(-1)^{pc(c)}$ from Theorem \ref{OHSachs1}, have the following sign pattern.
\begin{align*}
    \begin{array}{ccccc}
&+&&+& \\
+&&+&&- \\
&&+&& 
\end{array}
\end{align*}
From Theorem \ref{t:NonMonicIsZero} we know that summing over all edge-monic tail-equivalence classes produces the determinant of the Laplacian. Summing the contributors in this one class produces a value of $4$, and the determinant of the Laplacian would be calculated by repeating this process over all edge-monic tail-equivalence classes. Surprisingly, we show that the sum of a single edge-monic tail-equivalence class of an $n$-full oriented hypergraph is always equal to the determinant of the original $\{\pm{1}\}$-incidence-matrix, up to sign.

\subsection{Single Edge-monic Tail-equivalence Classes}
\label{ssec:OneEM}

Recall that $\mathfrak{A}_1$ is the set of edge-monic tail-equivalence classes, and let $\mathcal{A} \in \mathfrak{A}_1$ be any edge-monic tail-equivalence class (regardless of identifier). 

\begin{theorem}
    Let $G$ be an $n$-full oriented hypergraph. If $\mathcal{A} \in \mathfrak{A}_1$, then 
    \begin{align*}
    \left|\sum_{c\in\mathcal{A}}(-1)^{pc(c)}\right|=\lvert \det(\mathbf{H}_G) \rvert.
    \end{align*}
    \label{t:EMOnly}
\end{theorem}
\begin{proof}
We begin by converting the positive circle count $pc(c)$ into a negative circle count as follows:
\begin{align*}
    \dsum_{c\in\mathcal{A}}(-1)^{pc(c)} = \dsum_{c\in\mathcal{A}}(-1)^{tc(c) +nc(c)},
\end{align*}
where $nc(c)$ is the number of negative circles of $c$, and $tc(c)$ is the total number of circles in $c$.

Now, partition the total circle count into even circles $ec(c)$ and odd circles $oc(c)$ giving
\begin{align*}
    \dsum_{c\in\mathcal{A}}(-1)^{tc(c) +nc(c)} = \dsum_{c\in\mathcal{A}}(-1)^{ec(c)+oc(c)+nc(c)}.
\end{align*}
However, 
\begin{align*}
    \dsum_{c\in\mathcal{A}}(-1)^{ec(c)+oc(c)+nc(c)} = \dsum_{c\in\mathcal{A}}(-1)^{ec(c)} \dprod_{v \in V} \sigma(c(i_v))\sigma(c(j_v)).
\end{align*}

This last equality is easiest understood from right to left as it appears in \cite{OHSachs}. It is due to the fact every contributor corresponds to a permutation, and the value $\sigma(c(i_v))\sigma(c(j_v))$ is the adjacency/backstep entry of $\mathbf{L}_G$. Since $\mathbf{L}_G = \mathbf{D}_G - \mathbf{A}_G$ \cite{AH1} the parity equivalence is seen by factoring out a $-1$ for each adjacency (as it is subtracted) to get $(-1)^{oc(c)}$, and then factoring out a $-1$ for each negative adjacency to get $(-1)^{nc(c)}$. Observe that $(-1)^{ec(c)}$ is simply the sign of the permutation.

Taking the absolute value and separating the product gives
\begin{align*}
    \left|\sum_{c\in\mathcal{A}}(-1)^{pc(c)}\right| = \left| \dsum_{c\in\mathcal{A}}(-1)^{ec(c)} \dprod_{v \in V} \sigma(c(i_v)) \dprod_{v \in V}\sigma(c(j_v))\right|.
\end{align*}
However, since $\mathcal{A}$ is a tail-equivalence class, the tail product term
\begin{align*}
    \dprod_{v \in V} \sigma(c(i_v)) 
\end{align*}
is the same for each contributor and has magnitude $1$. Thus,
\begin{align*}
    \left|\sum_{c\in\mathcal{A}}(-1)^{pc(c)}\right| = \left| \dsum_{c\in\mathcal{A}}(-1)^{ec(c)}  \dprod_{v \in V}\sigma(c(j_v))\right|.
\end{align*}
This last term is $\det(\mathbf{H}_G)$ since every permutation is represented in an $n$-full tail-equivalence class, $(-1)^{ec(c)}$ is the sign of the corresponding permutation, and $\sigma(c(j_v))$ is the entry in the incidence matrix.   \qed
\end{proof}

Given an $n$-full oriented hypergraph $G$, let $\mathcal{A}^{+}_{\alpha}$ be the set of positive contributors in the edge-monic tail-equivalence class with identifier $\alpha$, and let $\mathcal{A}^{-}_{\alpha}$ be the set of negative contributors. Again, the sign of a contributor is determined by $(-1)^{pc(c)}$. Theorem \ref{t:EMOnly} gives the following immediate corollary.

\begin{corollary}
Let $G$ be an $n$-full oriented hypergraph, and $\mathcal{A}_{\alpha} \in \mathfrak{A}_1$, then 
    \begin{align*}
    \lvert \det(\mathbf{H}_G) \rvert = n! - 2\lvert  \mathcal{A}^{-}_{\alpha}\rvert = 2\lvert  \mathcal{A}^{+}_{\alpha}\rvert - n!,
    \end{align*}
    if the contributors $\mathcal{A}_{\alpha}$ sum to a non-negative integer, and 
    \begin{align*}
    \lvert \det(\mathbf{H}_G) \rvert = n! - 2\lvert  \mathcal{A}^{+}_{\alpha}\rvert = 2\lvert  \mathcal{A}^{-}_{\alpha}\rvert - n!,
    \end{align*}
    if the contributors $\mathcal{A}_{\alpha}$ sum to a non-positive integer.
\label{c:PminusNOneClass}
\end{corollary}
Turning our attention to all edge-monic tail-equivalence classes we can also examine the determinant of the Laplacian. Given an $n$-full oriented hypergraph $G$, let $\mathfrak{A}_{1}^{+}$ be the set of edge-monic tail-equivalence classes that sum to a positive integer and $\mathfrak{A}_{1}^{-}$ be the set of edge-monic tail-equivalence classes that sum to a negative integer.

\begin{lemma}
Let $G$ be an $n$-full oriented hypergraph corresponding to $\{\pm 1\}$-matrix $\mathbf{H}_G$ with non-zero determinant, then     
\begin{align*}
    \lvert \det(\mathbf{H}_G) \rvert = \lvert \mathfrak{A}_{1}^{+} \rvert - \lvert \mathfrak{A}_{1}^{-} \rvert = n! - 2\lvert \mathfrak{A}_{1}^{-} \rvert = 2\lvert \mathfrak{A}_{1}^{+} \rvert - n!.
    \end{align*}
\label{l:PminusNAllClass}
\end{lemma}
\begin{proof}
From Theorem \ref{t:NonMonicIsZero} we know
\begin{align*}
    \det(\mathbf{L}_{G}) &= \dsum\limits_{ \mathcal{A} \in \mathfrak{A}_1}\dsum\limits_{c \in \mathcal{A}} (-1)^{pc(c)} \\
    &= \dsum\limits_{ \mathcal{A} \in \mathfrak{A}_{1}^{+}}\dsum\limits_{c \in \mathcal{A}} (-1)^{pc(c)} + \dsum\limits_{ \mathcal{A} \in \mathfrak{A}_{1}^{-}}\dsum\limits_{c \in \mathcal{A}} (-1)^{pc(c)},
\end{align*}
From Theorem \ref{t:EMOnly} the inner sums become
\begin{align*}
    &= \dsum\limits_{ \mathcal{A} \in \mathfrak{A}_{1}^{+}}\lvert \det(\mathbf{H}_G) \rvert - \dsum\limits_{ \mathcal{A} \in \mathfrak{A}_{1}^{-}}\lvert \det(\mathbf{H}_G) \rvert \\
    &= \lvert \det(\mathbf{H}_G) \rvert \left(\lvert \mathfrak{A}_{1}^{+}  \rvert - \lvert \mathfrak{A}_{1}^{-} \rvert \right).
\end{align*}
The first equality is immediate from $\det(\mathbf{L}_{G}) = \det(\mathbf{H}_{G})^2$. While the last two equalities are immediate from the first equality and the symmetry of the edge-monic tail-equivalence classes of $n$-full oriented hypergraphs where $\lvert \mathfrak{A}_{1}^{+}  \rvert + \lvert \mathfrak{A}_{1}^{-} \rvert = n!$. \qed
\end{proof}

This underscores the symmetry of edge-monic tail-equivalence classes. The magnitude of the determinant of $\{\pm{1}\}$-matrices can be determined by the contributor signs in a single edge-monic tail-equivalence class, or by knowing the signs of the sums of each edge-monic tail-equivalence class individually. Adjacency-complements from Subsection \ref{ssec:EMIdentAndHeadTail} indicate that some signs are shared between the edge-monic tail-equivalence classes.

\section{Relationship to Maximum Determinants}
\label{sec:HadAndMax}

Hadamard's maximum determinant problem \cite{Had93} aims to find the maximum determinant of a matrix $\mathbf{H}$ of size $n$ with entries $+1$ and $-1$, and established a simple upper bound of $\lvert \det(\mathbf{H}) \rvert \leq n^{n/2}$. We refer the reader to \cite{Barba1933,browne2021survey, Ehlich1964,Paley1933OnOM,Wojtas1964} for classical approaches, bounds, and constructions, as well as a contemporary survey regarding the maximum determinant problem. As the oriented hypergraphic interpretation represents a clear deviation from these analyses, we conclude with a brief discussion on the oriented hypergraphic relationship to Hadamard's maximum determinant problem with the hope that the oriented hypergraphic interpretation provides new insights and approaches to complement existing methodologies. 

We make the simple observation in Subsection \ref{ssec:01Equiv} that the equivalence between the determinant of $\{\pm 1\}$-matrices and $\{0,+1\}$-matrices is related to the position of negative and positive digons that form a fundamental set of circles in the $n$-full oriented hypergraph. While Corollary \ref{c:PminusNOneClass} implies the maximal determinant is equivalent to optimizing the contributors of a specific sign in any edge-monic tail-equivalence class. Lemma \ref{l:PminusNAllClass} implies the maximal determinant is equivalent to minimizing the number of edge-monic tail-equivalence classes that sum to a negative integer. We conclude with Subsection \ref{ssec:ReclaimH} by reconstructing an incidence matrix given the signs of contributors in a single edge-monic tail-equivalence class. Given the extensive history of Hadamard's maximum determinant problem, further investigation is clearly needed in both relating these concepts to existing prior work, as well as understanding the highly symmetric combinatorial nature of $n$-full oriented hypergraphs and tail-equivalence classes.

\subsection{Fundamental Circles and Equivalence to $\{0,+1\}$-matrices}
\label{ssec:01Equiv}

Let $G = (V,E,I,\omega,\varsigma)$ be an incidence hypergraph. The \emph{cyclomatic number} of $G$ was introduced in \cite{OH1} and is
\begin{align*}
    \phi(G) = \lvert I \rvert - (\lvert V \rvert + \lvert E \rvert)+m,
\end{align*}
where $m$ is the number of connected components of $G$. The cyclomatic number of an $n$-full hypergraph is 
\begin{align*}
    \phi(G) = n^2 - (n + n)+1 = (n-1)^2.
\end{align*}
The incidence hypergraphic cyclomatic number is the minimal number of incidences whose deletion leaves a spanning circle-free incidence hypergraph, called a \emph{spanning incidence-forest}. Note, the deletion of all the incidences of an edge $e$ will result in $e$ becoming a $0$-edge, the edge-analog to an isolated vertex.

Paralleling our expectations from graph theory, the reintroduction of each incidence produces a unique circle in $G$ relative to spanning incidence-forest $F$ called a \emph{fundamental circle of $G$ with respect to $F$}. Again, while the resulting subhypergraph locally resembles a ``cycle'' or ``circuit'' of a graph we use ``circle'' as the resulting object need not be in the cycle space, nor need to be a matroid circuit --- this is an artifact of only having positive $2$-edges \cite{OH1}. 

Consider an $n$-full incidence hypergraph and remove the $(n-1)^2$ incidences between $v_k$ and $e_{\ell}$ for every $k,{\ell} \neq 1$ to produce a spanning incidence-tree. The reintroduction of a $(v_k,e_{\ell})$ incidence produces the digon that uses vertices $v_1$ and $v_k$ and edges $e_1$ and $e_{\ell}$. The set of these digons is called the \emph{fundamental bouquet} of $G$, and correspond to the minors of the form
\begin{align*}
    \left[ \begin{array}{cc}
        h_{1,1} & h_{1,\ell} \\ 
        h_{k,1} & h_{k,\ell}
        \end{array}
    \right]
\end{align*} 
in $\mathbf{H}_G$.

The fundamental bouquet provides a combinatorial interpretation on the known simple relationship between the determinants of $\{\pm 1\}$-matrices and the determinants of $\{0,1\}$-matrices as follows. Given a $\{\pm 1\}$-matrix $\mathbf{H}$ negate the necessary rows and columns so that the entries of the first row and first column are $+1$. Next, pivot on the $(1,1)$-entry so that the resulting $(1,1)$-minor only has entries $0$ or $-2$. Factoring out $-2$ from each row gives a matrix $\mathbf{H}'$ where
\begin{align*}
    \lvert \det(\mathbf{H}) \rvert = 2^{n-1} \lvert \det(\mathbf{H}')\rvert.
\end{align*}

\begin{example}
The following illustrates how to translate a $\{\pm 1\}$-matrix of size $4$ to a $\{0,1\}$-matrix while preserving the magnitude of the determinant.
\renewcommand{\arraystretch}{.6}
\[\begin{bmatrix}-1&1&-1&1\\-1&-1&-1&-1\\-1&1&1&-1\\1&1&-1&-1\end{bmatrix}\xrightarrow{negate}\begin{bmatrix}1&1&1&1\\1&-1&1&-1\\1&1&-1&-1\\1&-1&-1&1\end{bmatrix}\xrightarrow{pivot}\begin{bmatrix}1&1&1&1\\0&-2&0&-2\\0&0&-2&-2\\0&-2&-2&0\end{bmatrix}\]
\[\xrightarrow{minor}\begin{bmatrix}\textcolor{teal}{-2}&\textcolor{blue}{0}&\textcolor{red}{-2}\\0&-2&-2\\-2&-2&0\end{bmatrix}\xrightarrow{factor}\begin{bmatrix}1&0&1\\0&1&1\\1&1&0\end{bmatrix}.\]
\label{ex:reduction}
\end{example}

\begin{figure}[H]
    \centering
    \includegraphics{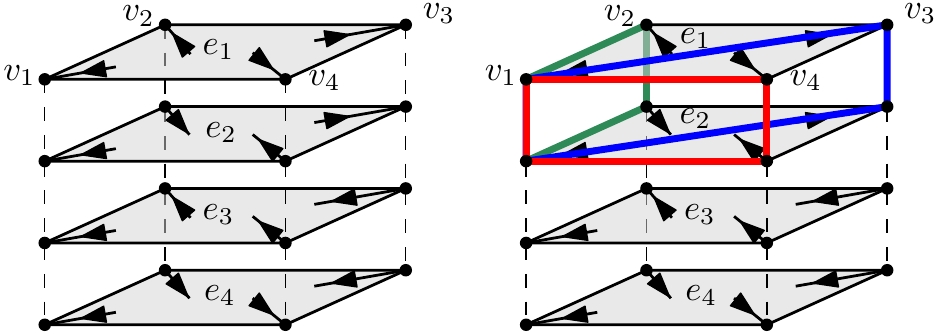}
    \caption{Left: A standardized $4$-full oriented hypergraph $G$; Right: The fundamental digons using edge $e_2$ correspond to entries in the minor.}
    \label{fig:Example1}
\end{figure}

An oriented hypergraph is said to be \emph{standardized} if all the incidences containing $v_1$ or $e_1$ are oriented $+1$. Figure \ref{fig:Example1} shows the oriented hypergraph for standardized matrix in Example \ref{ex:reduction} along with the fundamental bouquet digons using edge $e_2$. 

The entries in the $\{0,1\}$-matrix minor are determined by the signs of digons in the fundamental bouquet.

\begin{lemma}
    A digon in the fundamental bouquet of an $n$-full oriented hypergraph is positive (negative) if and only if the corresponding entry in the associated $\{0,+1\}$-matrix is $0$ ($+1$).
\end{lemma}
\begin{proof}
    Let $h_{k,\ell}$ be the $(k,\ell)$-entry of $\mathbf{H}$, and $h'_{k,\ell}$ be the $(k,\ell)$-entry of $\mathbf{H}'$, where $k,\ell \neq 1$. Observe that row and column negation always affects exactly $2$ incidences in any circle, hence it does not alter the signs of the digons in the fundamental bouquet, so we may assume the oriented hypergraph is standardized. 
    A digon in the fundamental bouquet is positive if and only if $h_{k,\ell} = 1$, and negative if and only if $h_{k,\ell} = -1$. The proof is complete with the observation that $h'_{k,\ell} = 0$ if and only if $h_{k,\ell} = 1$, and $h'_{k,\ell} = 1$ if and only if $h_{k,\ell} = -1$ \qed
\end{proof}

\subsection{Reconstructing $\{\pm 1\}$-matrices via Contributor Signs}
\label{ssec:ReclaimH}

In this subsection we assume all incidence matrices have been standardized so that the entries in the first row and column are $+1$, and from Theorem \ref{t:EMOnly} we assume we are working with $\mathcal{A}_{id}$. There are only $s(n,n-1) = \binom{n}{2}$ digons in a given edge-monic tail-equivalence class, so they cannot form a fundamental set of circles. However, the inclusion of the contributors corresponding to $3$-permutations of the form $(1 k \ell )$ with $1 < k < \ell$ produce another $s(n-1,n-2)$ circles, giving $s(n,n-1) + s(n-1,n-2) = (n-1)^2$ circles. We show the signs of the contributors corresponding to these permutations uniquely recovers the standardized incidence matrix. Hence, they must form a fundamental set of circles. 

Let $\Theta_k$ be the set of digon-contributors corresponding to permutation $(k \ell)$ with $1 \leq k < \ell$. These digons correspond to the $2 \times 2$ minors of the form
\begin{align*}
    \left[ \begin{array}{cc}
        h_{k,k} & h_{k,\ell} \\ 
        h_{\ell,k} & h_{\ell,\ell}
        \end{array}
    \right],
\end{align*}
as they iteratively travel down the identifier $id$. Alternate identifiers would produce appropriately adjusted minors.

\begin{theorem}
Given a standardized $n$-full oriented hypergraph $G$, the entries of the incidence matrix $\mathbf{H}_G$ are determined by the signs of the digon-contributors as well as those corresponding to $3$-permutations of the form $(1 k \ell )$ with $1 < k < \ell$ in $\mathcal{A}_{id}$ as follows:
\begin{enumerate}
    \item $h_{k,k} = - (-1)^{pc(c_{(1k)})}$,
    \item $h_{k,\ell} = (-1)^{pc(c_{(1 k \ell )})} \cdot (-1)^{pc(c_{(1k)})} \cdot (-1)^{pc(c_{(1\ell)})}$,
    \item $h_{\ell,k} = - (-1)^{pc(c_{(1 k \ell )})} \cdot (-1)^{pc(c_{(k \ell)})}$.
\end{enumerate}
\end{theorem}
\begin{proof}
Since $G$ is standardized all $h_{1,k} = h_{k,1} = +1$, and by Theorem \ref{t:EMOnly} assume we are working with $\mathcal{A}_{id}$. There is exactly one free entry (incidence) in each digon of $\Theta_1$ which corresponds to $h_{k,k}$, $k \geq 2$. The product of the incidences (entries) for these digons are 
\begin{align}
    h_{1,1} \cdot h_{1,k} \cdot h_{k,k} \cdot h_{k,1}  = h_{k,k}.
    \label{e:MainDiag}
\end{align}
The digons in the $\Theta_k$ ($1 < k < \ell$) correspond to the contributors $c_{(k \ell )}$, and have incidence (entry) product
\begin{align}
    h_{k,k} \cdot h_{k,\ell} \cdot h_{\ell,k} \cdot h_{\ell,\ell}.
    \label{e:BotDiag}
\end{align}
While the $3$-cycle in $c_{(1 k \ell )}$ has the following incidence product \begin{align}
    h_{1,1} \cdot h_{1,k} \cdot h_{k,k} \cdot h_{k,\ell} \cdot h_{\ell,\ell} \cdot h_{\ell,1} = h_{k,k} \cdot h_{k,\ell} \cdot h_{\ell,\ell}.
    \label{e:TopDiag}
\end{align}
Observe that the locally graphic nature of the $3$-cycle does not select the $3 \times 3$ minor, rather it selects only the $6$ entries of that minor that correspond to the incidences of the corresponding circle. 

Consider those contributors $c_{\pi}$ that correspond to the digons or single $3$-cycles of the form $(1 k \ell)$. Observe that if $c_{\pi}$ contains a single cycle of sign $\epsilon$, then $(-1)^{pc(c_{\pi})} = -\epsilon$. Thus, if $h_{k,k} = +1$ the digon in $c_{(1k)}$ is positive, and $(-1)^{pc(c_{(1k)})} = -1$; if $h_{k,k} = -1$ the digon in $c_{(1k)}$ is negative, and $(-1)^{pc(c_{(1k)})} = +1$; giving $(-1)^{pc(c_{(1k)})} = -h_{k,k}$. Thus, the entries on the main diagonal are determined by Equation \ref{e:MainDiag} and the digons from $\Theta_1$.

A similar argument applies to the entries above the main diagonal. Given the main diagonal entries, the only unknown entry in Equation \ref{e:TopDiag} corresponds to the $(k,\ell)$-entries. Using $(-1)^{pc(c_{(1k)})} = -h_{k,k}$ this reduces to $(-1)^{pc(c_{(1 k \ell )})} = h_{k,\ell} \cdot (-1)^{pc(c_{(1k)})} \cdot (-1)^{pc(c_{(1\ell)})}$. Finally, Equation \ref{e:BotDiag} contains the $(\ell , k)$-entries and uses the digons from $\Theta_k$, $k > 1$. A similar argument gives $(-1)^{pc(c_{(k \ell )})} = - (-1)^{pc(c_{(1k \ell)})} \cdot h_{\ell, k}$. Solving for the $h_{\ell, k}$ entry completes the proof.
\qed
\end{proof}

While the digons from $\Theta_1$ determine the signs of the main diagonal, it is surprising that the digons from the other $\Theta_k$ are necessarily used after the $3$-cycle contributors. While there seems to be a benefit for examining the contributors in this order, the ramifications are not clear. Also, choosing another identifier to initially populate a permutation in $\mathbf{H}$ other than the main diagonal would give rise to similar reconstructions, Theorem \ref{t:EMOnly} seems to indicate that it is unnecessary. In practice, we have demonstrated that for various $n$ you cannot simply require all the fundamental contributors to be a specific sign to produce a determinant that is maximal; the inherent contributor symmetries may play a role in this obstruction.

%%%%%%%%%%%%%%%%%%%%%%%%%%%%%%%%

\bibliographystyle{plain}
%\bibliography{references}

\begin{thebibliography}{10}

\bibitem{Barba1933}
Guido Barba.
\newblock Intorno al teorema di hadamard sui determinanti a valore massimo.
\newblock {\em Giorn. Mat. Battaglini}, 71:70--86, 1933.

\bibitem{browne2021survey}
Patrick Browne, Ronan Egan, Fintan Hegarty, and Padraig~O Cathain.
\newblock A survey of the hadamard maximal determinant problem, 2021.

\bibitem{Seth1}
Seth Chaiken.
\newblock A combinatorial proof of the all minors matrix tree theorem.
\newblock {\em SIAM Journal on Algebraic Discrete Methods 3}, 3:319--329, 1982.

\bibitem{OHSachs}
G.~Chen, V.~Liu, E.~Robinson, L.~J. Rusnak, and K.~Wang.
\newblock A characterization of oriented hypergraphic laplacian and adjacency
  matrix coefficients.
\newblock {\em Linear Algebra and its Applications}, 556:323 -- 341, 2018.

\bibitem{Reff6}
Luke Duttweiler and Nathan Reff.
\newblock Spectra of cycle and path families of oriented hypergraphs.
\newblock {\em Linear Algebra Appl.}, 578:251--271, 2019.

\bibitem{MR0267898}
Jack Edmonds and Ellis~L. Johnson.
\newblock Matching: {A} well-solved class of integer linear programs.
\newblock In {\em Combinatorial {S}tructures and their {A}pplications ({P}roc.
  {C}algary {I}nternat., {C}algary, {A}lta., 1969)}, pages 89--92. Gordon and
  Breach, New York, 1970.

\bibitem{Ehlich1964}
Hartmut Ehlich.
\newblock Determinantenabschätzungen für binäre matrizen.
\newblock {\em Math. Z.}, 83:123--132, 1964.

\bibitem{Had93}
J.~Hadamard.
\newblock {R\'{e}solution d'une question relative aux d\'{e}terminants}.
\newblock {\em Bulletin des Sciences Math\'{e}matiques.}, 17:240--246, 1893.

\bibitem{Mulas4}
Jürgen Jost and Raffaella Mulas.
\newblock Hypergraph laplace operators for chemical reaction networks.
\newblock {\em Advances in Mathematics}, 351:870 -- 896, 2019.

\bibitem{Reff7}
Ouail Kitouni and Nathan Reff.
\newblock Lower bounds for the laplacian spectral radius of an oriented
  hypergraph.
\newblock {\em Australasian Journal of Combinatorics}, 74(3):408--422, 2019.

\bibitem{Mulas3}
R.~Mulas.
\newblock Sharp bounds for the largest eigenvalue of the normalized hypergraph
  laplace operator.
\newblock {\em Mathematica Notes, to appear}, 2020.

\bibitem{Mulas1}
R.~Mulas.
\newblock Spectral classes of hypergraphs.
\newblock {\em ArXiv:2007.04273 [math.CO]}, 2020.

\bibitem{Mulas2}
R.~Mulas and D.~Zhang.
\newblock Spectral theory of laplace operators on chemical hypergraphs.
\newblock {\em ArXiv:2004.14671 [math.CO]}, 2020.

\bibitem{Paley1933OnOM}
R.~Paley.
\newblock On orthogonal matrices.
\newblock {\em Journal of Mathematics and Physics}, 12:311--320, 1933.

\bibitem{AH1}
N.~Reff and L.J. Rusnak.
\newblock An oriented hypergraphic approach to algebraic graph theory.
\newblock {\em Linear Algebra and its Applications}, 437(9):2262--2270, 2012.

\bibitem{Reff2}
Nathan Reff.
\newblock Spectral properties of oriented hypergraphs.
\newblock {\em Electron. J. Linear Algebra}, 27:373--391, 2014.

\bibitem{Reff3}
Nathan Reff.
\newblock Intersection graphs of oriented hypergraphs and their matrices.
\newblock {\em Australas. J. Combin.}, 65:108--123, 2016.

\bibitem{Reff5}
Nathan Reff and Howard Skogman.
\newblock A connection between {H}adamard matrices, oriented hypergraphs and
  signed graphs.
\newblock {\em Linear Algebra Appl.}, 529:115--125, 2017.

\bibitem{OHMTT}
L.~J. Rusnak, E.~Robinson, M.~Schmidt, and P.~Shroff.
\newblock Oriented hypergraphic matrix-tree type theorems and bidirected minors
  via boolean ideals.
\newblock {\em J Algebr Comb}, 49(4):461–--473, 2019.

\bibitem{OH1}
L.J. Rusnak.
\newblock Oriented hypergraphs: Introduction and balance.
\newblock {\em Electronic J. Combinatorics}, 20(3)(\#P48), 2013.

\bibitem{Wojtas1964}
W.~Wojtas.
\newblock On hadamard's inequality for the determinants of order non-divisible
  by 4.
\newblock {\em Colloq. Math}, 12:73--83, 1964.

\bibitem{OSG}
Thomas Zaslavsky.
\newblock Orientation of signed graphs.
\newblock {\em European J. Combin.}, 12(4):361--375, 1991.

\end{thebibliography}

\end{document}